\documentclass[graybox]{svmult}

\usepackage{type1cm}        
\usepackage{makeidx}         
\usepackage{graphicx}        
\usepackage{multicol}        
\usepackage[bottom]{footmisc}

\usepackage{newtxtext}       %
\usepackage[varvw]{newtxmath}       

\usepackage[numbers]{natbib}
\usepackage{multirow}

\makeindex             

\begin{document}
	
  \title*{Parallel iterative solvers for discretized reduced optimality systems}
	\titlerunning {Solvers for discretized reduced optimality systems}
	\author{Ulrich Langer, Richard L\"oscher, Olaf Steinbach, Huidong Yang}
        \institute{
          Ulrich Langer \at Institute of Numerical Mathematics, 
          JKU Linz,
          Austria, \email{ulanger@numa.uni-linz.ac.at} 
          \and 
          Richard L\"oscher \at Institut f\"{u}r Angewandte Mathematik, 
          TU Graz,
          Austria, \email{loescher@math.tugraz.at} 
          \and 
          Olaf Steinbach \at Institut f\"{u}r Angewandte Mathematik, TU Graz, Austria, \email{loescher@math.tugraz.at} 
          \and 
          Huidong Yang \at Faculty of Mathematics, University of Vienna, Austria,
          \email{huidong.yang@univie.ac.at} 
          }
	
	%
	%
	\maketitle

\abstract{
We propose, analyze, and test new iterative solvers 
for large-scale systems of linear algebraic equations arising from 
the 
finite element discretization of reduced optimality systems
defining the finite element approximations to the solution of
elliptic
tracking-type distributed optimal control problems 
with both the standard $L_2$  and the more general energy regularizations.
If we aim at 
an
approximation of the given desired state $y_d$ 
by the computed finite element state $y_h$   
that asymptotically differs from  $y_d$  in the order of the best $L_2$ approximation
under acceptable
costs for the control,
then the optimal choice of the regularization parameter $\varrho$ is 
linked
to the mesh-size $h$ by the relations
$\varrho=h^4$ and $\varrho=h^2$ for the $L_2$ and the energy regularization, respectively.
For this setting, we can construct efficient parallel iterative solvers for the reduced finite element optimality systems.
These results can be generalized to 
variable regularization parameters 
adapted to the local behavior 
of the mesh-size that can heavily change in case of adaptive 
mesh refinement. 
Similar results can be obtained for the space-time finite element discretization 
of the corresponding parabolic and hyperbolic optimal control problems.}

\section{Introduction}
\label{Sec:Introduction}

Let us first consider an abstract tracking-type, distributed Optimal Control Problem (OCP)
of the form: Find the state  $y_\varrho \in Y$ and the control $u_\varrho \in U$  
minimizing the cost functional
\begin{equation}
\label{Eqn:AbstractCostFunctional1}
J_\varrho(y_\varrho,u_\varrho):= \frac{1}{2} \|y_\varrho - y_d\|_H^2 + \frac{\varrho}{2} \| u_\varrho \|_U^2
                              = \frac{1}{2} \|y_\varrho - y_d\|_H^2 + \frac{1}{2} \| \sqrt{\varrho} \, u_\varrho \|_U^2
\end{equation}
subject to (s.t.) the state equation
\begin{equation}
\label{Eqn:AbstractStateEquation1}
B y_\varrho = u_\varrho 
\;\; \mbox{in} \;\; U \subseteq P^*,
\end{equation}
where $y_d \in H$ denotes the given desired state (target), $\varrho > 0$ is a suitably chosen 
regularization parameter that also 
affects
the energy cost $\| u_\varrho \|_U^2$ for the control $u_\varrho$ as source term in \eqref{Eqn:AbstractStateEquation1}, and
$X=Y, P,  U, H$ are Hilbert spaces equipped with the corresponding norms $\|\cdot\|_X$
and scalar products $(\cdot,\cdot)_X$.
We assume that $Y \subset H \subset Y^*$ and $P \subset H \subset P^*$ form Gelfand triples of 
Hilbert spaces, and that $B: Y \rightarrow P^*$ is an isomorphism,
where $X^*$ denotes the dual space of $X$ 
with the duality product 
$\langle \cdot,\cdot \rangle: X^* \times X \rightarrow \mathbb{R}$ 
that is nothing but  the extension of the scalar product $(\cdot,\cdot)_H$
in $H$ for $X=Y$ and $X=P$.
We are interested in the cases $U=H$ and $U=P^*$. 
Optimal control problems of the form \eqref{Eqn:AbstractCostFunctional1}-\eqref{Eqn:AbstractStateEquation1}
with many applications were already investigated in 
the classical monograph \cite{DD28:Lions:1968a} by Lions and
the more recent book \cite{DD28:Troeltzsch:2010a} by Tr\"oltzsch,
where additional constraints of the form $u_\varrho \in U_{\text{ad}} \subset U$ imposed on the control 
$u_\varrho$ are  permitted. The unique solvability of such kind of OCPs is based on
the unique solvability of the state equation, i.e. $y_\varrho = B^{-1} u_\varrho$,
the strong convexity of the quadratic cost functional and the assumption that the admissible set $U_{\text{ad}}$
is a non-empty, convex, and closed subset of $U$;
cf. Theorem~2.16 in \cite{DD28:Troeltzsch:2010a}.
Here we only consider the case $U_{\text{ad}}=U$. Then the unique solution 
$(y_\varrho,u_\varrho) \in Y \times U$ of the 
OCP \eqref{Eqn:AbstractCostFunctional1} - \eqref{Eqn:AbstractStateEquation1} 
can also be extracted from
the unique solution $(y_\varrho,p_\varrho,u_\varrho) \in Y \times P \times U$ 
of the first-order optimality system (OS)
\begin{equation}
\label{Eqn:OS}
B y_\varrho = u_\varrho,\, B^* p_\varrho = y_\varrho - y_d,\, p_\varrho + A_{1/\varrho}^{-1}u_\varrho = 0, 
\end{equation}
where the self-adjoint and elliptic regularization operator 
$A_{1/\varrho}: P \rightarrow P^*$ is defined by the regularization via 
the Riesz representation of the control. For $U=P^*$, we have 
$A_{1/\varrho} = \varrho^{-1} A$ and  $\|u\|_{P^*}^2 = \langle u,A^{-1}u\rangle$,
whereas $A=I$ (canonical embedding operator) for $U=H$. 
Formally, we will write 
$\| \sqrt{\varrho} \, u \|_U^2 = \langle u,A_{1/\varrho}^{-1}u\rangle$
that allows us to permit variable $\varrho$.
Eliminating $u = - A_{1/\varrho}p$ from 
\eqref{Eqn:OS}, 
we arrive at the equivalent reduced OS (ROS) written as saddle point problem: 
Find $(y_\varrho,p_\varrho) \in Y \times P$ such that 
        \begin{equation}
        \label{Eqn:ROS}
        \begin{bmatrix}
            A_{1/\varrho} & B\\
            B^*             & -I
        \end{bmatrix}
        \begin{bmatrix}
            p_\varrho\\
            y_\varrho
        \end{bmatrix}
        =
        \begin{bmatrix}
            0\\
            -y_d
        \end{bmatrix}
        \quad \mbox{in}\;P^* \times Y^*.
        \end{equation}
%
%
%

Typical examples are elliptic, parabolic, and hyperbolic OCPs, 
where the state equation \eqref{Eqn:AbstractStateEquation1} is given by 
an elliptic boundary value problem (BVP), 
a parabolic initial boundary value problem (IBVP), 
and a hyperbolic IBVP,
respectively. 
In this paper, we will focus on the parallel numerical solution of elliptic OCPs,
where the Dirichlet boundary value problem for the Poisson equation 
serves us as model problem for the state equation. 
%
We are primarily interested in efficient parallel solvers for algebraic systems 
arising from the finite element (FE) discretization of the ROS \eqref{Eqn:ROS} 
or the corresponding primal Schur complement 
when the FE discretization and the regularization 
are balanced in  an asymptotically optimal way.

\section{Elliptic Optimal Control Problems}
\label{Sec:EOCP}

As mentioned 
above,
we focus on elliptic OCPs of the form 
\eqref{Eqn:AbstractCostFunctional1}-\eqref{Eqn:AbstractStateEquation1}
defined by the following specifications:
$Y = P = H_0^1(\Omega)$, 
$Y^* = P^* = H^{-1}(\Omega)$, $H=L_2(\Omega)$, and 
$B = - \Delta: Y = H_0^1(\Omega) \rightarrow P^* = H^{-1}(\Omega)$ is defined by the variational identity 
\begin{equation}
\label{Eqn:EOCP:B}
\langle B y , p \rangle = (\nabla y , \nabla p)_{L_2(\Omega)},   
\;\; \forall y \in Y = H_0^1(\Omega),\; \forall p \in P = H_0^1(\Omega),
\end{equation}
where $\Omega \subset \mathbb{R}^d,\, d=1,2,3,$ denotes the $d$-dimensional
computational domain that is supposed to be bounded and Lipschitz.  
Throughout the paper, we use the usual notations for Lebesgue and Sobolev
spaces. In \cite{DD28:NeumuellerSteinbach:2021a},
Neum\"uller and Steinbach showed that $\|y_\varrho - y_d\|_{L_2(\Omega)}$ behaves 
like $O(\varrho^{s/r})$ provided that
$y_d \in H_0^s(\Omega) := [L_2(\Omega),H_0^1(\Omega)]_s$
for some $s \in [0,1]$, where $r=4$ ($U=L_2(\Omega)$),
and $r=2$ ($U=H^{-1}(\Omega)$), respectively.

The FE 
Galerkin 
discretization of 
the ROS 
\eqref{Eqn:ROS} reads as follows:
Find $y_h=y_{\varrho h} \in Y_h = S_h^1(\mathcal{T}_h)\cap Y =\text{span}\{\varphi_{hk}\}_{k=1}^{n_h} \subset Y$
     and 
     $p_h=p_{\varrho h} \in P_h = S_h^1(\mathcal{T}_h)\cap P = \text{span}\{\psi_{hi}\}_{i=1}^{m_h} \subset P$
such that 
%
\begin{equation}
 \label{Eqn:EOCP:FEROS}
 \langle A_{1/\varrho}p_h,q_h\rangle + \langle B y_h,q_h\rangle = 0
 \;\mbox{and}\;
 \langle B^* p_h,v_h\rangle - 
  (y_h,v_h)_{L_2(\Omega)} = - (y_d,v_h)_{L_2(\Omega)}
\end{equation}
for all $q_h \in P_h$ and $ v_h \in Y_h$,
where $S_h^1(\mathcal{T}_h)$ is nothing but the continuous, piecewise linear 
FE space defined on some shape-regular triangulation $\mathcal{T}_h$ of $\Omega$.
Here $Y_h=P_h$,  $n_h = m_h$ and $\varphi_{hk} = \psi_{hk}$ since $Y=P=H_0^1(\Omega)$.
Once the basis is chosen, the FE ROS leads to the symmetric and indefinite (SID) algebraic system
        \begin{equation}
        \label{Eqn:EOCP:FEROS-GS}
        \begin{bmatrix}
            A_{1/\varrho, h} & B_h\\
            B^T_h             & -M_h
        \end{bmatrix}
        \begin{bmatrix}
            \mathbf{p}_h\\
            \mathbf{y}_h
        \end{bmatrix}
        =
        \begin{bmatrix}
            \mathbf{0}_h\\
            -\mathbf{y}_{dh}
        \end{bmatrix}
        \end{equation}
that can further be  reduced to the symmetric and positive definite (SPD) Schur complement (SC) system
        \begin{equation}
        \label{Eqn:EOCP:FEROS-SC}
        (B^T_h A_{1/\varrho, h}^{-1} B_h + M_h) \mathbf{y}_h = \mathbf{y}_{dh}, 
        \end{equation}
where $\mathbf{y}_h= [y_k]_{k=1}^{n_h} \in \mathbb{R}^{n_h}$ and 
      $\mathbf{p}_h= [y_k]_{k=1}^{n_h} \in \mathbb{R}^{n_h}$
are the nodal FE vectors corresponding to the FE functions 
$y_h \in Y_h = P_h$ and $p_h \in P_h = Y_h$ 
via the FE isomorphism $\mathbf{y}_h, \mathbf{p}_h\leftrightarrow  y_h, p_h$.
The SPD SC system \eqref{Eqn:EOCP:FEROS-SC} can be solved by the 
PCG method 
provided that some good preconditioner for the SC is available, and 
the application of 
$A_{\varrho h}^{-1}$ to some vector 
performs
in asymptotically optimal complexity.


Let us first briefly 
review some results concerning the
standard $L_2$ regularization where $U=H=L_2(\Omega)$
yielding $A_{1/\varrho, h} = \varrho^{-1} M_h$.
Then we can prove the estimate
\begin{equation}
 \label{Eqn:EOCP:DiscretizationErrorEstimate}
 \|y_{\varrho h} - y_d\|_{L_2(\Omega)} \le c h^s \|y_d\|_{H^s(\Omega)}
\end{equation}
provided that $\varrho = h^4$ \cite{DD28:LangerLoescherSteinbachYang:2023CMAM}.
Estimate \eqref{Eqn:EOCP:DiscretizationErrorEstimate} remains true 
for $A_{1/\varrho, h} = \varrho^{-1} D_h$, where $D_h = \mbox{lump}(M_h)$ denotes the lumped mass matrix 
\cite{DD28:LangerLoescherSteinbachYang:2023arXiv:2304.14664}.
In both cases, the lumped mass matrix $D_h$ is spectrally equivalent to the Schur complement
$S_h = B^T_h A_{1/\varrho, h}^{-1} B_h + M_h$ 
provided that $\varrho = h^4$. More precisely,
the spectral equivalence inequalities
\begin{equation}
  \label{Eqn:EOCP:SpectralEquivalenceInequalities}
  (d+2)^{-1} D_h \le M_h \le S_h := B^T_h A_{1/\varrho, h}^{-1} B_h + M_h \le (c_\text{\tiny inv}^r + 1) M_h \le (c_\text{\tiny inv}^r + 1) D_h
\end{equation}
hold with $r=4$, where the positive constant $c_\text{\tiny inv}$ is defined by the inverse inequality
%
%
$\|\nabla v_h\|_{L_2(\Omega} \le c_\text{\tiny inv} h^{-1}\|v_h\|_{L_2(\Omega)}$ for all $v_h \in Y_h$.
We refer to \cite{DD28:LangerLoescherSteinbachYang:2023arXiv:2304.14664}
for the proof of 
\eqref{Eqn:EOCP:SpectralEquivalenceInequalities}.
In particular, for targets $y_d$ not belonging to $Y$, e.g. discontinuous targets, 
it may be useful to use adaptivity and variable, mesh-adapted regularizations
of the form 
%
\begin{equation}
 \label{Eqn:EOCP:DiffusionRegularization} 
 \varrho(x) = h_\tau^2, \forall x \in \tau, \forall \tau \in \mathcal{T}_h.
\end{equation}
%
Now $A_{1/\varrho, h} = M_{1/\varrho}$, where $M _{1/\varrho}$ is defined by
$(M_{1/\varrho}\mathbf{p}_h,\mathbf{q}_h) = ((1/\varrho) p_h,q_h)_{L_2(\Omega)}$,
see \cite{DD28:LangerLoescherSteinbachYang:2023:BerichteTUGraz}  for details.

In this contribution, we will
focus on the more non-standard $H^{-1}$ regularization where $U=P^*=H^{-1}(\Omega)$
yielding $A_{1/\varrho, h} = K_{1/\varrho, h}$. 
The regularization matrix $K_{1/\varrho, h}$ is now a diffusion stiffness matrix defined by 
$(K_{1/\varrho,h}\mathbf{p}_h,\mathbf{q}_h) = ((1/\varrho) \nabla p_h,\nabla q_h)_{L_2(\Omega)}$,
with the mesh-adapted choice \eqref{Eqn:EOCP:DiffusionRegularization} 
%
%
of the regularization $\varrho$.
This kind of diffusion regularization
was completely analyze in 
\cite{DD28:LangerLoescherSteinbachYang:2022arXiv:2209.08811}.
In particular, the spectral equivalence inequalities \eqref{Eqn:EOCP:SpectralEquivalenceInequalities}
are valid with $r=2$. 
Thus, in principle, the SPD SC system \eqref{Eqn:EOCP:FEROS-SC} can be solved by the
PCG with the diagonal matrix $D_h$ as preconditioner. 
However, in contrast to the $L_2$ regularization, the matrix $A_{1/\varrho, h} = K_{1/\varrho, h}$
cannot be replace by a diagonal matrix in order to obtain a fast matrix-by-vector multiplication 
without disturbing the accuracy of the discretization. 
One loophole would be the use of inner iteration to invert $K_{1/\varrho, h}$,
e.g. (algebraic) multigrid iteration as we used in some numerical experiments in 
\cite{DD28:LangerLoescherSteinbachYang:2022arXiv:2209.08811}. 
In order to avoid inner iterations down to the discretization error,
we can solve the larger SID system \eqref{Eqn:EOCP:FEROS-GS} 
by some SID solver like MINRES, Bramble-Pasciak PCG, or GMRES.
Then we only need a preconditioner for $K_{1/\varrho, h}$ and $S_h$.
In the case of a constant regularization parameter $\varrho = h^2$,
we surprisingly observe that $A_{1/\varrho, h} = K_{1/\varrho, h} = (1/\varrho) K_h = (1/\varrho) B_h$,
and, therefore, $S_h = B^T_h A_{1/\varrho, h}^{-1} B_h + M_h = \varrho K_h + M_h$,
and the SC system \eqref{Eqn:EOCP:FEROS-SC} simplifies to a diffusion problem
of the form $(\varrho K_h + M_h)\mathbf{y}_h = \mathbf{y}_{dh}$.
This system can easily be solved by PCG with the diagonal preconditioner $D_h$.
%
       
Inspired by this observation, we now propose to solve the diffusion equation 
\begin{equation}
\label{Eqn:EOCP:DiffusionEquation}
        (\varrho \nabla y, \nabla v)_{L_2(\Omega)} + (y,v)_{L_2(\Omega)} = (y_d,v)_{L_2(\Omega)}\; \forall v \in Y = H^1_0(\Omega)
\end{equation}
when we want to choose variable regularization parameters $\varrho=\varrho(x)$ 
in connection with adaptive FE discretization.
Obviously, the diffusion equation  \eqref{Eqn:EOCP:DiffusionEquation} is nothing but 
the first-order OS for minimizing the cost functional
\begin{equation*}
{\widetilde J}(y) = 0.5 [\|y-y_d\|^2+   \|\sqrt{\varrho} \nabla y\|^2]
= 0.5 [\|y-y_d\|^2 +   \langle B^{-*}A_{\varrho}B^{-1}u,u\rangle] = {\widetilde J}(y,u)
\end{equation*}
instead of the original cost functional
\begin{equation*}
 J(y) = 0.5 [\|y-y_d\|^2 +   \langle B^* A_{1/\varrho}^{-1}By,y\rangle]]
                          = 0.5 [\|y-y_d\|^2 +   \langle A_{1/\varrho}^{-1}u,u\rangle] = J(y,u),                         
\end{equation*}
where $B^{-*} := (B^{-1})^*$, and the subscript $L_2(\Omega)$ is here omitted from the norms.
%
We note that $\langle A_{1/\varrho}^{-1}u,u\rangle \le \langle B^{-*}A_{\varrho}B^{-1}u,u\rangle$
with ``$=$'' instead of ``$\le$'' for constant $\varrho$.

The FE discretization of \eqref{Eqn:EOCP:DiffusionEquation} now leads to the SPD system
\begin{equation}
 \label{Eqn:EOCP:(K+M)y=rhs}
 (K_{\varrho h} + M_h) \mathbf{y}_h = \mathbf{y}_{dh},
\end{equation}
where the SPD diffusion stiffness matrix $K_{\varrho h}$ is defined by the identity
\begin{equation}
 \label{Eqn:EOCP:K}
 (K_{\varrho h}\mathbf{y}_h, \mathbf{v}_h) =  (\varrho \nabla y_h, \nabla v_h)_{L_2(\Omega)}, \;
 \forall \, \mathbf{y}_h, \mathbf{v}_h\leftrightarrow  y_h, v_h \in Y_h
\end{equation}
For the diffusion regularization \eqref{Eqn:EOCP:DiffusionRegularization},
the discretization error $\|y_h - y_d\|_{L_2(\Omega)}$
can be analyzed in the similar way as was done for the original approach 
in \cite{DD28:LangerLoescherSteinbachYang:2022arXiv:2209.08811}.
Moreover, the SPD system \eqref{Eqn:EOCP:(K+M)y=rhs} is much more simpler 
that the original systems \eqref{Eqn:EOCP:FEROS-GS} or \eqref{Eqn:EOCP:FEROS-SC},
and can easily be solved by PCG. Indeed, the diffusion regularization \eqref{Eqn:EOCP:DiffusionRegularization}
ensures that the system matrix $K_{\varrho h} + M_h$ of \eqref{Eqn:EOCP:(K+M)y=rhs} 
fulfills the same spectral equivalence inequalities \eqref{Eqn:EOCP:SpectralEquivalenceInequalities}
as the SC $B^T_h A_{1/\varrho, h}^{-1} B_h + M_h$ with the same $r=2$.
Thus, we can solve \eqref{Eqn:EOCP:(K+M)y=rhs} by means of the PCG preconditioned 
by the diagonal matrix $D_h = \mbox{lump}(M_h)$, or by another diagonal approximation of 
the mass matrix $M_h$ like $\mbox{diag}(M_h)$. This leads to an asymptotically optimal 
solver for a fixed relative accuracy. This solver can be used in a nested iteration setting 
on a sequence of refined meshes $\mathcal{T}_\ell$, $\ell = 1,2,\ldots,L$, 
in such a way that, at each level $\ell$, we compute final iterates that differs
from the target $y_d$ in the order of the discretization error with respect to 
the $L_2$ norm in optimal complexity. The the whole solution process has optimal complexity.
We note that, due to \eqref{Eqn:EOCP:SpectralEquivalenceInequalities}, 
the PCG converges in the $L_2$ norm, and the number of nested iterations is de facto 
constant over the discretization levels. It is clear that both the non-nested and 
the nested PCG solver can easily be parallelized.
These convergence properties and parallel performance will impressively be confirmed by our numerical 
experiments presented in the next section.




\section{Numerical experiments}
\label{Sec:NUMA}


We here consider the discontinuous desired state
\begin{equation}
 \label{Eqn:NUMA:Desired State}
 y_d = 1\;\mbox{in}\; (0.25,0.75)^3,\;\mbox{and}\; 
 y_d = 0\;\mbox{in}\; \overline{\Omega} \setminus (0.25,0.75)^3,
\end{equation}
where $\Omega = (0,1)^3 \subset \mathbb{R}^{d=3}$.
We note that 
this discontinuous desired state 
$y_d \in L_2(\Omega)$ does not belong to the 
state space $H_0^1(\Omega)$, but $y_d \in H^{1/2 - \varepsilon}(\Omega)$
for any $\varepsilon > 0$.

We start with a uniform decomposition of the
domain $\Omega=(0,1)^3$ into $24,576$
tetrahedral elements $\tau \in \mathcal{T}_1$. 
The coarsest mesh  $\mathcal{T}_1$ contains  $4,913$ 
vertices, 
and the mesh size $h=h_1=0.0625$.
A sequence of  meshes $\mathcal{T}_\ell$ at the levels $\ell=2,...,6$ is
generated by successive uniform refinement. 
On the finest level $\ell=L=6$, there are $135,005,697$ vertices, $h=1.953125$e$-3$,
and $\varrho=h^2\approx 3.81$e$-6$. 
In the adaptive algorithm, we have
chosen the locally adapted $\varrho_\tau=h_\tau^2$ on each $\tau \in \mathcal{T}_h$.
The adaptivity is driven by the localization of the error
$\|\tilde{y}_\ell-y_d\|_{L_2(\Omega)}$.
System \eqref{Eqn:EOCP:(K+M)y=rhs} is solved
by PCG with the preconditioner $\mbox{diag}(M_h)$
The PCG iteration stops when the
relative preconditioned residual is reduced by a factor of $10^6$.
The parallel implementation is based on the opensource
MFEM 
({\tt https://mfem.org/}),
and is running on the HPC cluster
RADON1 
({\tt https://www.oeaw.ac.at/ricam/hpc}).

In the nested iteration, 
we run the PCG until the relative preconditioned residual reaches $10^{-6}$ at the coarsest level $\ell=1$. 
At the refined levels $\ell=2,3,...$, an adaptive tolerance
$\alpha[n_\ell/n_{\ell-1}]^{-\beta/3}$, $\ell=2,3,...$, is adopted for controlling
the relative preconditioned residual, with $n_\ell$ being the number of degrees
of freedom (\#Dofs) at the mesh level $\ell$. Here $\alpha$ is a scaling
factor, e.g., $0.5$ and $0.25$ for the uniform and adaptive refinement, respectively. 
The parameter $\beta$ is directly related to the  convergence rate of the 
discretization error $\|y_\ell - y_d\|_{L_2(\Omega)}$. 
For our example, $\beta$ is $0.5$ and $0.75$ for the uniform and the adaptive refinement, 
respectively.

\begin{figure}[h]
  \centering
    \includegraphics[width=1\textwidth]{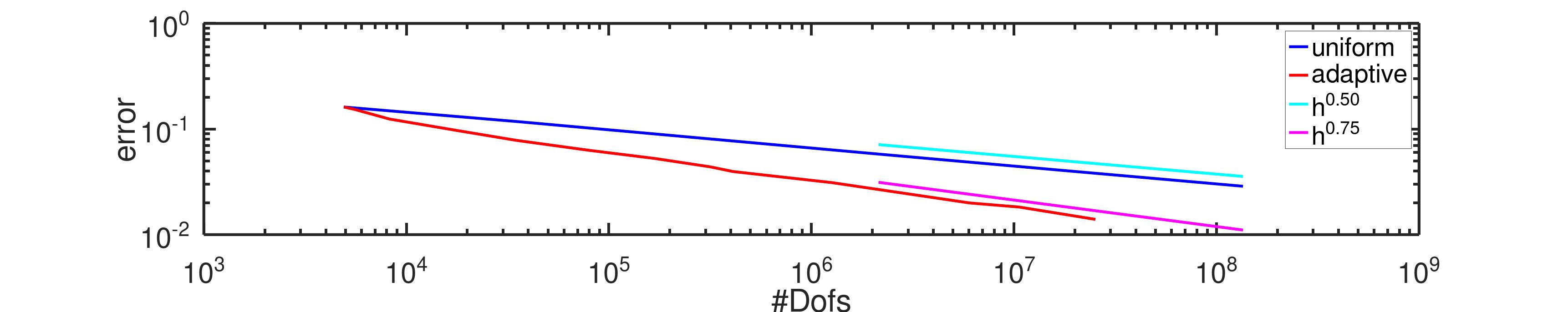}
    \caption{Non-nested iteration: Convergence history for uniform and adaptive refinements.
      }
    \label{fig:errorcomp_nonnested}
\end{figure}

The \#Dofs, the error $\|y_d - y_\ell^{k_\ell}\|_{L_2(\Omega)}$,
the number Its of PCG iterations, and the
corresponding computational time (Time) in seconds at each level $\ell$ are illustrated
in Table~\ref{tab:solver_nonnested} for non-nested iterations 
starting with the initial guess $y_\ell^{0} = 0$.
All computations were performed using $512$ cores.
It is easy to observe the expected convergence
rate $0.5$ for the uniform refinement and $0.75$ for the adaptive
refinement; see also Fig.~\ref{fig:errorcomp_nonnested}.
The robustness of our preconditioner is confirmed by the
constant PCG iteration numbers for both uniform and adaptive refinements.
The efficiency of the parallelization is here demonstrated by very 
small
computational costs at all levels. 
The adaptive method outperforms the uniform one in terms of computational time needed for
achieving a similar accuracy. 

\begin{figure}[h]
  \centering
    \includegraphics[width=1\textwidth]{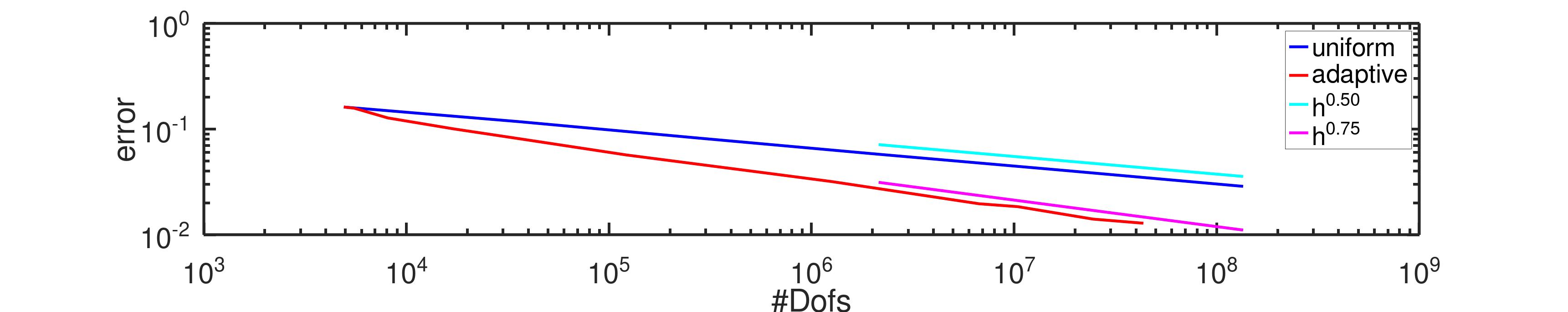}
    \caption{Nested iteration: Convergence history for uniform and adaptive refinements.
      }
    \label{fig:errorcomp_nested}
\end{figure}

For the nested iteration where the initial guess is interpolated 
from the coarser mesh, the same behavior can be observed 
from Table~\ref{tab:solver_nested} which presents the same data as Table~\ref{tab:solver_nonnested}.
The number of PCG iterations as well as the computational time
are reduced by a factor of $8$ approximately without loss of accuracy; 
cf. Fig.~\ref{fig:errorcomp_nested}. 

Finally, in  Table \ref{tab:parallel_solver_nonnested} for the non-nested
iterations 
and Table \ref{tab:parallel_solver_nested} for nested iterations,
we compare computational time in seconds with respect to both the number of cores and
number of refinement levels on uniform refinement.
We clearly see the strong (rowwise) and weak (diagonal) scalability of our parallel solvers.

\begin{table}[h]
  {\small
    \begin{tabular}{|l|r|l|r|r|l|r|r|}
      \hline
      \multirow{2}{*}{$\ell$}&\multicolumn{3}{c|}{Adaptive} &
      \multicolumn{4}{c|}{Uniform} \\  \cline{2-8}
      &\#Dofs& error &Its (Time)
      &\#Dofs& error &eoc & Its (Time)\\
      \hline
      $1$&$4,913$&$1.61$e$-1$&$20$ ($6.0$e$-3$ s)&$4,913$&$1.61$e$-1$& $-$ &$20$ ($6.0$e$-3$ s)\\
      $2$&$5,532$&$1.54$e$-1$&$24$ ($8.3$e$-3$ s)&$35,937$&$1.17$e$-1$ &$0.46$ & $23$ ($7.4$e$-3$ s)\\
      $3$&$8,255$&$1.24$e$-1$&$25$ ($8.5$e$-3$ s)&$274,625$&$8.26$e$-2$&$0.51$ & $23$ ($9.1$e$-3$ s)\\
      $4$&$18,013$&$9.65$e$-2$&$26$ ($9.9$e$-3$ s)&$2,146,689$&$5.79$e$-2$&$0.51$ & $22$ ($2.0$e$-2$ s)\\
      $5$&$35,055$&$7.80$e$-2$&$26$ ($1.1$e$-2$ s)&$16,974,593$&$4.07$e$-2$&$0.51$ & $22$ ($2.0$e$-1$ s)\\
      $6$&$80,381$&$6.28$e$-2$&$27$ ($1.3$e$-2$ s)&$135,005,697$&$2.87$e$-2$&$0.50$ & $22$ ($1.4$e$-0$ s)\\
      $7$&$167,982$&$5.27$e$-2$&$27$ ($1.4$e$-2$ s)&&&& \\
      $8$&$316,839$&$4.48$e$-2$&$27$ ($1.7$e$-2$ s)&&&& \\
      $9$&$410,144$&$3.96$e$-2$&$27$ ($2.0$e$-2$ s)&&&&\\
      $10$&$1,264,336$&$3.11$e$-2$&$28$ ($2.8$e$-2$ s)&&&&\\
      $11$&$6,043,649$&$2.00$e$-2$&$27$ ($1.3$e$-1$ s)&&&&\\
      $12$&$10,590,586$&$1.83$e$-2$&$28$ ($2.6$e$-1$ s)&&&&\\
      $13$&$25,217,222$&$1.40$e$-2$&$27$ ($5.0$e$-1$ s)&&&&\\
      \hline
    \end{tabular}
    \caption{Non-nested iteration: Error, number Its of PCG iterations and computational time (Time) 
      for adaptive and uniform refinements, where eoc denotes the experimental order of convergence.
    } 
    \label{tab:solver_nonnested}
    }
\end{table}

\begin{table}[h]
  {\small
    \begin{tabular}{|l|r|l|r|r|l|r|r|}
      \hline
      \multirow{2}{*}{$\ell$}&\multicolumn{3}{c|}{Adaptive} &
      \multicolumn{4}{c|}{Uniform} \\  \cline{2-8}
      &\#Dofs& error &Its (Time)
      &\#Dofs& error &eoc & Its (Time)\\
      \hline
      $1$&$4,913$&$1.61$e$-1$&$20$ ($6.0$e$-3$ s)&$4,913$&$1.61$e$-1$& $-$ &$20$ ($6.0$e$-3$ s)\\
      $2$&$5,532$&$1.58$e$-1$&$2$ ($9.4$e$-4$ s)&$35,937$&$1.18$e$-1$&$0.46$ & $2$ ($8.9$e$-4$ s)\\
      $3$&$8,096$&$1.27$e$-1$&$3$ ($1.2$e$-3$ s)&$274,625$&$8.22$e$-2$&$0.52$& $3$ ($1.4$e$-3$ s)\\
      $4$&$17,166$&$1.00$e$-1$&$3$ ($1.2$e$-3$ s)&$2,146,689$&$5.77$e$-2$&$0.52$ & $3$ ($3.7$e$-3$ s)\\
      $5$&$34,134$&$9.23$e$-2$&$2$ ($1.4$e$-3$ s)&$16,974,593$&$4.09$e$-2$&$0.50$ & $2$ ($3.0$e$-2$ s)\\
      $6$&$73,530$&$6.59$e$-2$&$3$ ($2.1$e$-3$ s)&$135,005,697$&$2.87$e$-2$&$0.52$ & $2$ ($1.8$e$-1$ s)\\
      $7$&$121,987$&$5.68$e$-2$&$2$ ($1.8$e$-3$ s)&&&& \\
      $8$&$624,691$&$3.80$e$-2$&$2$ ($2.5$e$-3$ s)&&&& \\
      $9$&$1,260,214$&$3.18$e$-2$&$2$ ($3.0$e$-3$ s)&&&&\\
      $10$&$6,719,190$&$1.96$e$-2$&$2$ ($1.6$e$-2$ s)&&&&\\
      $11$&$10,426,031$&$1.84$e$-2$&$2$ ($2.8$e$-2$ s)&&&&\\
      $12$&$24,509,387$&$1.40$e$-2$&$1$ ($3.4$e$-2$ s)&&&&\\
      $13$&$43,437,311$&$1.28$e$-2$&$3$ ($1.3$e$-1$ s)&&&&\\
      \hline
    \end{tabular}
    \caption{Nested iteration: Error, number Its of PCG iterations and computational time (Time) 
      for adaptive ($\alpha=0.25$, $\beta=0.75$) and uniform refinements ($\alpha=0.5$, $\beta=0.5$). 
      %
      } 
    \label{tab:solver_nested}
    }
\end{table}

  \begin{table}[h]
    {\small
      \begin{tabular}{|c|c|c|c|c|c|c|c|}
        \hline\multirow{2}{*}{$\ell$}&\multicolumn{6}{c|}{\#Cores}\\  \cline{2-7}
        &16&32&64&128&256&512\\
        \hline
        $2$&$23$ ($9.3$e$-3$ s)&-&-&-&-&-\\ 
        $3$&$23$ ($8.7$e$-2$ s)&$23$ ($4.4$e$-2$ s)&$23$ ($1.8$e$-2$ s)&$23$ ($1.2$e$-2$ s)&-&-\\
        $4$&$22$ ($6.1$e$-1$ s)&$22$ ($3.2$e$-1$ s)&$22$ ($1.8$e$-1$ s)&$22$ ($9.7$e$-2$ s)&$22$ ($4.9$e$-2$ s)&$22$ ($2.0$e$-2$ s)\\
        $5$&$22$ ($4.9$e$-0$ s)&$22$ ($2.6$e$-0$ s)&$22$ ($1.3$e$-0$ s)&$22$ ($6.9$e$-1$ s)&$22$ ($3.7$e$-1$ s)&$22$ ($2.0$e$-1$ s)\\
        $6$&-&-&- &$22$ ($5.5$e$-0$ s)&$22$ ($2.7$e$-0$ s)&$22$ ($1.4$e$-0$ s) \\ \hline
      \end{tabular}
      \caption{Number of non-nested PCG iterations and computational time 
      for uniform refinement.}
      \label{tab:parallel_solver_nonnested}
    }
  \end{table}

  \begin{table}[h]
    {\small
      \begin{tabular}{|c|c|c|c|c|c|c|c|}
        \hline\multirow{2}{*}{$\ell$}&\multicolumn{6}{c|}{\#Cores}\\  \cline{2-7}
        &16&32&64&128&256&512\\
        \hline
        $2$&$2$ ($1.1$e$-3$ s)&-&-&-&-&-\\ 
        $3$&$3$ ($1.4$e$-2$ s)&$3$ ($7.1$e$-3$ s)&$3$ ($2.8$e$-3$ s)&$3$ ($1.6$e$-3$ s)&-&-\\
        $4$&$3$ ($1.0$e$-1$ s)&$3$ ($5.6$e$-2$ s)&$3$ ($2.9$e$-2$ s)&$3$ ($1.7$e$-2$ s)&$3$ ($8.9$e$-3$ s)&$3$ ($3.7$e$-3$ s)\\
        $5$&$2$ ($6.2$e$-1$ s)&$2$ ($3.3$e$-1$ s)&$2$ ($1.6$e$-1$ s)&$2$ ($8.8$e$-2$ s)&$2$ ($4.7$e$-2$ s)&$2$ ($3.0$e$-2$ s)\\
        $6$&-&-&- &$2$ ($6.9$e$-1$ s)&$2$ ($3.4$e$-1$ s)&$2$ ($1.8$e$-1$ s) \\ \hline
      \end{tabular}
      \caption{Number of nested PCG iterations and computational time 
               for uniform refinement.}
      \label{tab:parallel_solver_nested}
    }
  \end{table}

\section{Conclusions, Generalizations, and Outlook}\label{Sec:???}

We have presented a new diffusion regularization that leads to simple 
diffusion equation as first order optimality condition. 
The corresponding FE equations can very efficiently be solved 
by parallel PCG with diagonal preconditioners provided the 
the regularization is appropriately chosen. 
It is possible to include box constraints for both the control 
and the state \cite{DD28:GanglLoescherSteinbach:2023arXiv:2306.15316}.
These techniques can be extended to parabolic and hyperbolic OCPs.

\begin{acknowledgement}
  We thank the RICAM for providing us with the high
  performance cluster Radon1. The last author is supported
  by the Christian Doppler Research Association. 
\end{acknowledgement}		

\bibliographystyle{plain} 
\bibliography{LangerLoescherSteinbachYang_miniXY}
\end{document}